\documentclass{birkmult}
\usepackage{graphicx,amssymb} 
\usepackage{epsfig}
\usepackage{mathrsfs}

\usepackage{bbm}

\title[Fractal zeta functions and complex dimensions]{Fractal zeta functions and complex dimensions: 
A general higher-dimensional theory
}

\author[M.\ L.\ Lapidus]{Michel L.\ Lapidus}
\address{University of California, Riverside\\ 
Department of Mathematics\\
California 92521-0135\\ 
USA}
\email{lapidus@math.ucr.edu}

\thanks{The research of Michel L.~Lapidus was partially supported by the National Science
Foundation under grants ~DMS-0707524 and DMS-1107750, as well as by the Institut des Hautes \' Etudes Scientifiques (IH\' ES) where the first author was a visiting professor in the Spring of 2012 while part of this research was completed. This paper is based on a plenary lecture given by the first author (MLL) at the International Conference on ``Fractal Geometry and Stochastics'' held in Tabarz, Germany, in March 2014. Goran Radunovi\'c and Darko \v{Z}ubrini\'c express their gratitude to the Ministry of Science of the Republic of Croatia for its support.}

\author[G.\ Radunovi\'c]{Goran Radunovi\'c}
\address{University of Zagreb, Faculty of Electrical Engineering and Computing\\ 
Department of Applied Mathematics\\
 Unska 3, 10000 Zagreb\\ 
Croatia}
\email{goran.radunovic@fer.hr} 

\author[D.\ \v Zubrini\'c]{Darko \v Zubrini\'c}
\address{University of Zagreb, Faculty of Electrical Engineering and Computing\\ 
Department of Applied Mathematics\\  
Unska 3, 10000 Zagreb\\ 
Croatia}
\email{darko.zubrinic@fer.hr}

\begin{document}

\begin{abstract}
In 2009, the first author introduced a class of zeta functions, called `distance zeta functions', which has enabled us to extend the existing theory of zeta functions of fractal strings and sprays (initiated by the first author and his collaborators in the early 1990s) to arbitrary bounded (fractal) sets in Euclidean spaces of any dimensions.
A closely related tool is the class of `tube zeta functions', defined using the tube function 
of a fractal set. These zeta functions exhibit deep connections with Minkowski contents 
and upper box (or Minkowski) dimensions, as well as, more generally, with the complex dimensions of fractal sets.
In particular, the abscissa of (Lebesgue, i.e., absolute) convergence of the distance zeta function coincides with the upper box dimension of a set.
We also introduce a class of transcendentally quasiperiodic sets, and describe their construction based on a sequence of carefully chosen 
generalized Cantor sets with two auxilliary parameters. 
As a result, we obtain a family of  ``maximally hyperfractal'' compact sets and relative fractal drums (i.e., such that the associated fractal zeta functions have 
a singularity at every point of the critical line of convergence).
Finally, we discuss the general fractal tube formulas and the Minkowski measurability criterion obtained by the authors in the context of relative fractal drums (and, in particular, of bounded subsets of $\mathbb{R}^N$).
\end{abstract}


\date{}

\maketitle


\newtheorem{theorem}{Theorem}[section]
\newtheorem{cor}[theorem]{Corollary}
\newtheorem{prop}[theorem]{Proposition}
\newtheorem{lemma}[theorem]{Lemma}

\theoremstyle{remark}
\newtheorem{remark}[theorem]{Remark}
\newtheorem{defn}[theorem]{Definition}
\newtheorem{example}[theorem]{Example}
\newtheorem{problem}[theorem]{Problem}
\newtheorem{exercise}[theorem]{Exercise}

\newcommand{\re}{\operatorname{Re}}
\newcommand{\im}{\operatorname{Im}}

\font\csc=cmcsc10

\def\esssup{\mathop{\rm ess\,sup}}
\def\essinf{\mathop{\rm ess\,inf}}
\def\wo#1#2#3{W^{#1,#2}_0(#3)}
\def\w#1#2#3{W^{#1,#2}(#3)}
\def\wloc#1#2#3{W_{\scriptstyle loc}^{#1,#2}(#3)}
\def\osc{\mathop{\rm osc}}
\def\var{\mathop{\rm Var}}
\def\supp{\mathop{\rm supp}}
\def\Cap{{\rm Cap}}
\def\norma#1#2{\|#1\|_{#2}}

\def\C{\Gamma}

\let\text=\mbox

\catcode`\@=11
\let\ced=\c
\def\a{\alpha}
\def\b{\beta}
\def\c{\gamma}
\def\d{\delta}
\def\g{\lambda}
\def\o{\omega}
\def\q{\quad}
\def\n{\nabla}
\def\s{\sigma}
\def\div{\mathop{\rm div}}
\def\sing{{\rm Sing}\,}
\def\singg{{\rm Sing}_\ty\,}

\def\A{{\cal A}}
\def\F{{\cal F}}
\def\H{{\cal H}}
\def\W{{\bf W}}
\def\M{{\cal M}}
\def\N{{\cal N}}

\def\eR{{\bf R}}
\def\eN{{\bf N}}
\def\Ze{{\bf Z}}
\def\Qe{{\bf Q}}
\def\Ce{{\bf C}}

\def\ty{\infty}
\def\e{\varepsilon}
\def\f{\varphi}
\def\:{{\penalty10000\hbox{\kern1mm\rm:\kern1mm}\penalty10000}}
\def\ov#1{\overline{#1}}
\def\DD{\Delta}
\def\O{\Omega}
\def\pa{\partial}

\def\st{\subset}
\def\stq{\subseteq}
\def\pd#1#2{\frac{\pa#1}{\pa#2}}
\def\sgn{{\rm sgn}\,}
\def\sp#1#2{\langle#1,#2\rangle}

\newcount\br@j
\br@j=0
\def\q{\quad}
\def\gg #1#2{\hat G_{#1}#2(x)}
\def\inty{\int_0^{\ty}}
\def\od#1#2{\frac{d#1}{d#2}}

\def\bg{\begin}
\def\eq{equation}
\def\bgeq{\bg{\eq}}
\def\endeq{\end{\eq}}
\def\bgeqnn{\bg{eqnarray*}}
\def\endeqnn{\end{eqnarray*}}
\def\bgeqn{\bg{eqnarray}}
\def\endeqn{\end{eqnarray}}

\def\bgeqq#1#2{\bgeqn\label{#1} #2\left\{\begin{array}{ll}}
\def\endeqq{\end{array}\right.\endeqn}

\def\abstract{\bgroup\leftskip=2\parindent\rightskip=2\parindent
        \noindent{\bf Abstract.\enspace}}
\def\endabstract{\par\egroup}

\def\udesno#1{\unskip\nobreak\hfil\penalty50\hskip1em\hbox{}
             \nobreak\hfil{#1\unskip\ignorespaces}
                 \parfillskip=\z@ \finalhyphendemerits=\z@\par
                 \parfillskip=0pt plus 1fil}
\catcode`\@=11

\def\cal{\mathcal}
\def\eR{\mathbb{R}}
\def\eN{\mathbb{N}}
\def\Ze{\mathbb{Z}}
\def\Qu{\mathbb{Q}}
\def\Ce{\mathbb{C}}

\def\osd{\mathrm{osd}\,}

\def\sdim{\mbox{\rm s-dim}\,}
\def\sd{\mbox{\rm sd}\,}

\def\res{\operatorname{res}}

\def\L{{\cal L}}
\def\po{\mathcal{P}}
\def\I{\mathbbm{i}}
\def\E{\mathrm{e}}
\def\D{\mathrm{d}}
\def\qs{\q}
\def\avM{\widetilde{\mathcal{M}}}




\section{Introduction}\label{intro}

This article provides a short survey of some of the recent advances in the theory of fractal zeta functions and the associated higher-dimensional theory of complex dimensions, valid for arbitrary bounded subsets of Euclidean spaces and developed in the forthcoming research monograph \cite{fzf}, entitled {\em Fractal Zeta Functions and Fractal Drums$:$ Higher-Dimensional Theory of Fractal Dimensions}. 
(See also the research articles [LapRa\v Zu2--6] and the survey article [LapRa\v Zu7].)

\smallskip

The theory of zeta functions of fractal strings, initiated by the first author 
in the early 1990s and described in an extensive research monograph \cite{lapidusfrank12}, joint with M.\ van Frankenhuijsen (see also the references therein),
was given an unexpected impetus in
 2009, when a new class of zeta functions, called `distance zeta functions', was discovered (also by the first author).\footnote{For fractal string theory and the associated one-dimensional theory of complex dimensions, as well as for the extensions to higher-dimensional fractal sprays (in the sense of [LapPo3]), we refer the reader to the research monographs [Lap-vFr1--3] along, for example, with the articles [Lap1--6, LapPo1--3, LapMa1--2, Fal1, HeLap, HamLap, Tep1--2, LapRo, LapL\'eRo, L\'eMen, LapLu, LapLu-vFr1--2, LapPe1--2, LapPeWi1--2, LalLap1--2, LapRo\v Zu, RaWi, ElLapMacRo, HerLap1--2]. We refer, in particular, to [Lap-vFr3, {\S}12.2.1 and Chap.\ 13] for a survey of some of the recent developments of the theory, prior to [LapRa\v Zu1--7].}
Since distance zeta functions are associated with arbitrary bounded (fractal) sets in Euclidean spaces of any dimension (see Definition \ref{defn}),
they clearly represent a valuable tool connecting the geometry of fractal sets with complex analysis.
This interplay is described in~[LapRa\v Zu1--7], where the foundations of the theory of fractal zeta functions have been laid.
In this paper, by `fractal zeta functions' we mean the following three classes of zeta functions: zeta functions of fractal strings (and, more generally, of fractal sprays), distance zeta functions and tube zeta functions of bounded subsets of $\eR^N$, with $N\ge1$, although some other classes may appear as well, like zeta functions of relative fractal drums in $\eR^N$ and spectral zeta functions; see {\S}\ref{fzfrfd} below and [LapRa\v Zu1, Chap.\ 4].
The theory of fractal zeta functions exhibits very interesting connections with the Minkowski contents and dimensions of fractal sets; see Theorems \ref{pole1} and~\ref{pole1mink_tilde}.

\smallskip

Like fractal string theory, which the present theory of fractal zeta functions extends to arbitrary dimensions (as well as to ``relative fractal drums'' in $\eR^N$), the work described here should eventually have applications to various aspects of harmonic analysis, fractal geometry, dynamical systems, geometric measure theory and analysis on nonsmooth spaces, number theory and arithmetic geometry, mathematical physics and, more speculatively, to aspects of condensed matter physics and cosmology. Some of the more mathematical applications of the theory are described in \cite{fzf}, as well as in [LapRa\v Zu2--7], but a variety of potential applications remain to be explored or even imagined.

The basic property of the distance zeta function of a fractal set, described in Theorem \ref{an}, is that its abscissa of (absolute or Lebesgue) convergence is equal to the upper box dimension $D$ of the set. 
Under some mild hypotheses, $D$ is always a singularity; see part~$(b)$ of Theorem~\ref{an}.
Furthermore, assuming that $D$
is a pole, then it is simple.
Moreover, the residue of the distance zeta function computed at $D$ is, up to a multiplicative constant, between the corresponding upper and lower Minkowski contents.
A similar statement holds for the tube zeta function.
(See Theorems~\ref{pole1} and~\ref{pole1mink_tilde}, respectively.)

\smallskip

In addition, according to part $(b)$ of Theorem~\ref{an},
 under some mild assumptions on a bounded set $A$, the abscissa of (Lebesgue, i.e., {\em absolute}) convergence of its distance zeta function coincides not only with $D$, but also with the abscissa of {\em holomorphic continuation} of the zeta function. 
\smallskip

We stress that if $D:=\ov\dim_BA<N$, all the results concerning the distance zeta functions have exact counterparts for the tube zeta functions, and vice versa.
In other words, the fractal zeta functions introduced in [LapRa\v Zu1--7] contain essentially the same information. In practice, however, it is often the case that one of the fractal zeta functions is better suited for the given situation under consideration.
\smallskip 

In \S\ref{merom_ext}, we discuss the existence and the construction of a suitable meromorphic continuation of the distance (or tube) zeta function of a fractal set, both in the Minkowski measurable case (Theorem \ref{measurable}) and a frequently encountered instance of Minkowski nonmeasurable case (Theorem \ref{nonmeasurable}). We will illustrate the latter situation by computing the fractal zeta function and the associated complex dimensions of the Sierpi\'nski carpet; see Proposition \ref{sierpinski_carpet0} and Example \ref{ex2} when $N=2$ or $3$, respectively. Many other examples are provided in \cite{fzf} and [LapRa\v Zu2--6], where are calculated, in particular, the complex dimensions of the higher-dimensional Sierpi\'nski gaskets and carpets in $\eR^N$, for any $N\ge2$. 

\smallskip

In {\S}\ref{qp_sets}, we introduce the so-called transcendentally $n$-quasiperiodic sets, for any integer $n\ge2$ (that is, roughly speaking, the sets possessing $n$ quasiperiods; see Definition \ref{quasiperiodic}), and describe the construction of $2$-quasiperiodic sets,  based on carefully chosen generalized Cantor sets with two parameters, introduced in Definition \ref{Cma}; see Theorem \ref{trans}. It is also possible to construct $n$-quasiperiodic sets, for any $n\ge2$, and even $\ty$-quasiperiodic sets, that is, sets which possess infinitely many quasiperiods; see \S\ref{maxhf} below and  \cite[{\S}4.6]{fzf}. 

\smallskip

In {\S}\ref{fzfrfd}, we introduce the notion of a relative fractal drum $(A,\O)$ (which represents a natural extension of the notion of bounded fractal string and of bounded set). We also introduce the corresponding relative distance and tube zeta functions $\zeta_A(\,\cdot\,,\O)$ and $\widetilde\zeta_A(\,\cdot\,,\O)$, and study their properties. It is noteworthy that the relative box dimension $\dim_B(A,\O)$ can be naturally defined as a real number, which may also assume negative values, including $-\ty$. 

\smallskip

In {\S}\ref{mmc}, we address the question of reconstructing the tube function $t\mapsto|A_t\cap\O|$ of a relative fractal drum $(A,\O)$, and thereby of obtaining a general ``fractal tube formula'' expressed in terms 
of the complex dimensions of $(A,\O)$ (defined as the poles of a suitable meromorphic extension of the relative distance zeta function $\zeta_A(\,\cdot\,,\O)$).
The corresponding tube formulas are obtained in \cite[Chap.\ 5]{fzf} and [LapRa\v Zu6] (announced in [LapRa\v Zu5]), as well as illustrated by a variety of examples. The example of the three-dimensional Sierpi\'nski carpet is given in Example \ref{ex2}. Moreover, towards the end of {\S}\ref{mmc}, we explain how to deduce from our general tube formulas (and significantly extend) earlier results obtained for fractal strings (in [Lap-vFr1--3]) and, especially, for fractal sprays and self-similar tilings (in [LapPe2] and [LapPeWi1]).

\medskip

In closing this introduction, we recall some basic notation and terminology which will be needed in the sequel. First of all, in order to avoid trivial special cases, we assume implicitly that all bounded subsets of $\eR^N$ under consideration in the statements of the theorems are nonempty.
Assume that $A$ is a given bounded subset of $\eR^N$ and let $r$ be a fixed real number. We define the {\em upper} and {\em lower $r$-dimensional Minkowski contents} of $A$, respectively, by
$$
\M^{*r}(A):=\limsup_{t\to0^+}\frac{|A_t|}{t^{N-r}},\q \M_*^r(A):=\liminf_{t\to0^+}\frac{|A_t|}{t^{N-r}},
$$
where $A_t$ denotes the Euclidean $t$-neighborhood of $A$ (namely, $A_t:=\{x\in\eR^N:d(x,A)<t\}$) and $|A_t|$ is the $N$-dimensional Lebesgue measure of $A_t$.
The {\em upper} and {\em lower box $(\mathrm{or}$ Minkowski$)$ dimensions} of $A$ are then defined, respectively, by
$$
\ov\dim_BA:=\inf\{r\in\eR:\M^{*r}(A)=0\},\q \underline\dim_BA:=\inf\{r\in\eR:\M_*^r(A)=0\}.
$$
It is easy to check that $0\leq\underline{\dim}_BA\leq\ov{\dim}_BA\leq N$.
Furthermore, if $A$ is such that $\underline\dim_BA=\ov\dim_BA$, then this common value is denoted by $\dim_BA$ and is called the {\em box $(\mathrm{or}$ Minkowski$)$ dimension} of $A$.
Moreover, if $A$ is such that, for some $D\in[0,N]$, we have $0<\M_*^D(A)\le\M^{*D}(A)<\ty$ (in particular, then $\dim_BA$ exists and $D=\dim_BA$), we say that {\em $A$ is  Minkowski nondegenerate}. If $\M_*^D(A)=\M^{*D}(A)$, then this common value is denoted by $\M^D(A)$ and called the {\em Minkowski content} of $A$. Finally, assuming that $A$ is such that $\M^D(A)$ exists and $0<\M^D(A)<\ty$, we say that {\em $A$ is Minkowski measurable}.\footnote{We note that the notion of Minkowski dimension was introduced (for noninteger values) by Bouligand \cite{Bou} in the late 1920s (without making a clear distinction between the lower and upper limits), while the notions of (lower and upper) Minkowski content, Minkowski measurability and Minkowski nondegeneracy were introduced, respectively, in~\cite{federer}, \cite{stacho} and \cite{rae}. (See also [Lap1, LapPo1] and, especially, [LapPo2--3], along with \cite{lapidusfrank12}, for the latter notions.) For general references on the notion of Minkowski (or box) dimension (from different points of view), we refer, for example, to [Fed2], [Fal2], [Tri], [Mat] and [Lap-vFr3].}

Throughout this paper, given $\alpha\in\eR\cup\{\pm\ty\}$, we denote by $\{\re s>\alpha\}$ the corresponding open right half-plane in the complex plane, defined by $\{s\in\Ce\,:\,\re s>\alpha\}$. (In particular, if $\a=\pm\ty$, $\{\re s>\a\}$ is equal to $\emptyset$ or $\Ce$, respectively.)
Similarly, given any $\alpha\in\eR$, we denote by $\{\re s=\alpha\}$ the corresponding vertical line $\{s\in\Ce\,:\,\re s=\alpha\}$.

\section{Distance and tube zeta functions}\label{disttube}

Let us introduce a new class of zeta functions, defined by the first author in 2009, which
extends the notion of geometric zeta functions of bounded fractal strings to bounded subsets of Euclidean spaces of arbitrary dimensions.

\smallskip

\begin{defn}[{[LapRa\v Zu1,2]}]\label{defn}  Let $A$ be a bounded subset of $\eR^N$ and let $\delta$ be a fixed positive real number.
Then, the {\em distance zeta function} $\zeta_A$ of $A$ is defined by
\begin{equation}\label{z}
\zeta_A(s):=\int_{A_\delta}d(x,A)^{s-N}\D x,
\end{equation}
for all $s\in\Ce$ with $\re s$ sufficiently large. Here, $d(x,A):=\inf\{|x-y|:y\in A\}$ denotes the usual Euclidean distance from $x$ to $A$.
Furthermore, the integral is taken in the sense of Lebesgue, and hence, is absolutely convergent.\footnote{For simplicity, we implicitly assume throughout this paper that $|A|=0$; the case when $|A|>0$ is discussed in \cite{fzf}.}
\end{defn}

\smallskip

\begin{remark}
Since the difference of any two distance zeta functions of the same set $A$ corresponding to two different values of $\delta$ is an entire function,\footnote{This is an easy consequence of the fact that $d(x,A)\in[\d_1,\d_2]$ for all $x\in A_{\d_2}\setminus A_{\d_1}$ with $0<\d_1<\d_2<\ty$.}  
it follows that the dependence of the distance zeta function $\zeta_A$ on $\delta>0$ is inessential, in the sense that the poles (of meromorphic extensions) of $\zeta_A$, as well as their multiplicities, do not depend on the choice of $\delta$.
\end{remark} 

The key for understanding the behavior of the distance zeta function $\zeta_A$ consists in understanding the Lebesgue integrability
of the function $A_\d\ni x\mapsto d(x,A)^{\re s-N}$, where $s\in \Ce$ is fixed.\footnote{Indeed, note that $|d(x,A)^{s-N}|=d(x,A)^{\re s-N}$ for all $x\in A_\d$.} (We shall soon see   that $\re s$ should be sufficiently large.) More precisely, we are interested in the Lebesgue integrability of the function $x\mapsto d(x,A)^{-\c}$ defined on $A_\d$, where $\gamma:=N-\re s$ and $s$ is a fixed complex number. Since the function is clearly bounded (and hence, integrable) for $\c\le0$, it suffices to consider the case when $\c>0$, that is, when $\re s<N$.

Let us recall a useful and little known result due to Harvey and Polking, stated implicitly on page 42 of \cite{acta}, in which a sufficient condition for Lebesgue integrability is expressed in terms of the upper box dimension. {\em If $A$ is any nonempty bounded subset of $\eR^N$, then the following implication holds}:\footnote{Moreover, if we assume that $D:=\dim_BA$ exists, $D<N$ and $\M_*^D(A)>0$, then the converse implication holds as well; see \cite[Thm.\ 4.3]{rae}. (See also [\v Zu2, Thm.\ 4.1(b)].)}
\begin{equation}\label{iff0}
\c<N-\ov\dim_BA\q\Longrightarrow\q\int_{A_\d} d(x,A)^{-\c}\D x<\ty.
\end{equation}

\begin{remark}\label{princip}
The sufficient condition for the Lebesgue (i.e., absolute) integrability of the function $A_\d \ni x\mapsto d(x,A)^{s-N}$  
in the Harvey--Polking result in \eqref{iff0}, becomes
$\c:=N-\re s<N-\ov\dim_BA$, that is, $\re s>\ov\dim_BA$. In other words, $\zeta_A(s)$ is well defined for all $s\in\Ce$ 
in the open right half-plane $\{\re s>\ov\dim_BA\}$.
\end{remark}

\smallskip

The distance zeta function of a bounded set represents a natural extension of the notion of {\em geometric zeta function} $\zeta_{\mathcal{L}}$, associated with a bounded fractal string $\mathcal{L}=(\ell_j)_{j\ge1}$ (introduced by the first author and his collaborators\footnote{See, especially, [Lap2--3], [LapPo1--3], [LapMa1--2] and [HeLap].} in the early 1990s and extensively studied in [Lap-vFr1--3] 
and the relevant references therein):
\begin{equation}
\zeta_{\mathcal{L}}(s):=\sum_{j=1}^\ty(\ell_j)^s,
\end{equation}
for all $s\in\Ce$ with $\re s$ sufficiently large.
Here, a {\em bounded fractal string} $\mathcal{L}$ is defined as a nonincreasing infinite sequence of positive real numbers $(\ell_j)_{j\ge1}$ such that $\ell:=\sum_{j\ge1}\ell_j<\ty$.
Alternatively, $\mathcal L$ can be viewed as a bounded open subset $\O$ of $\eR$, in which case the $\ell_j$s are the lengths of the connected components (open intervals) of $\O$, written in nonincreasing order (so that $\ell_j\downarrow0$ as $j\to\ty$).

An important first result concerning $\zeta_{\mathcal L}$ (first observed in [Lap2--3], using a result from \cite{BesTay}) is that its abscissa of (absolute) convergence coincides with $\tilde D$ (the inner Minkowski dimension of $\mathcal L$ or, equivalently, of its fractal boundary $\pa\O$), defined by $\tilde D:=\ov\dim_B(\pa\O,\O)$; see definition \eqref{rbdim} below. For a direct proof of this statement, see [Lap-vFr3, Thm.\ 1.10] or [Lap-vFr3, Thm.\ 13.111] and [LapLu-vFr2]. In light of the next comment, it can be readily shown that part $(a)$ of Theorem \ref{an} below extends this result to arbitrary compact subsets of Euclidean spaces in any dimension; see [LapRa\v Zu1,2].

It is easy to see that the distance zeta function $\zeta_{A_{\mathcal{L}}}$ of the set
$$
A_{\mathcal{L}}:=\bigg\{a_k:=\sum_{j=k}^\ty\ell_j:k\ge1\bigg\}\stq[0,\ell],
$$
associated with $\mathcal L$, 
 and the geometric zeta function $\zeta_{\mathcal{L}}$ are connected by the following simple relation:
\begin{equation}\label{dist_geo}
\zeta_{A_{\mathcal{L}}}(s)=u(s)\,\zeta_{\mathcal{L}}(s)+v(s),
\end{equation}
for all complex numbers $s$ such that $\re s$ is sufficiently large, where $u$ and $v$ are holomorphic on $\Ce\setminus\{0\}$ and $u$ is nowhere vanishing. 
In particular, due to Theorem \ref{an} below, it follows that the abscissae of convergence  of the distance zeta function $\zeta_A$ and 
of the geometric zeta function $\zeta_{\mathcal{L}}$ coincide, and that the corresponding  poles located on the critical line
$\{\re s=\ov\dim_BA_{\mathcal L}\}$ (called {\em principal complex dimensions} of $\mathcal L$ or, equivalently, of $A_{\mathcal L}$), as well as their multiplicities, also coincide. 
The exact same results hold if $A_{\mathcal L}$ is replaced by $\pa\O$, the boundary of $\O$, where $\O$ is any geometric realization of $\mathcal L$ by a bounded open subset of $\eR$.
For more details, see [LapRa\v Zu1,2].

\smallskip

Before stating Theorem~\ref{an}, we need to introduce some terminology and notation, which will also be used in the remainder of the paper.

\smallskip

Given a meromorphic function (or, more generally, an arbitrary complex-valued function) $f=f(s)$, initially defined on some domain $U\subseteq\Ce$, we denote by $D_{\rm hol}(f)$ the unique extended real number (i.e., $D_{\rm hol}(f)\in\eR\cup\{\pm\ty\}$) such that $\{\re s>D_{\rm hol}(f)\}$ is the {\em maximal} open right half-plane (of the form $\{\re s>\alpha\}$, for some $\alpha\in\eR\cup\{\pm\ty\}$) to which the function $f$ can be holomorphically extended.\footnote{By using the principle of analytic continuation, it is easy to check that $D_{\rm hol}(f)$ and $\mathcal{H}(f)$ are well defined; see~\cite[{\S}2.1]{fzf}.}
This maximal (i.e., largest) half-plane is denoted by $\mathcal H(f)$ and called the {\em half-plane of holomorphic continuation} of $f$.

If, in addition, the function $f=f(s)$ is assumed to be given by a {\em tamed Dirichlet-type integral} (or DTI, in short),\footnote{This is the case of the classic (generalized) Dirichlet series and integrals~\cite{serre,Pos}, the classic arithmetic zeta functions (see, e.g., \cite[App.\ A]{lapidusfrank12} and [Lap 4, Apps.\ B, C \& E]), as well as of the geometric zeta functions of fractal strings studied in [Lap-vFr1--3] and of all the fractal zeta functions considered in this paper and in [LapRa\v Zu1--7].} of the form
\begin{equation}\label{41/4}
f(s):=\int_E\varphi(x)^s\D\mu(x),
\end{equation}
for all $s\in\Ce$ with $\re s$ sufficiently large, where $\mu$ is a (positive or complex) local (i.e., locally bounded) Borel measure on a given (measurable) space $E$ and 
\begin{equation}\label{tamed}
0\le\varphi(x)\le C\q \mbox{for\q$|\mu|$-a.e. $x\in E$,}
\end{equation}
where $C\ge0$,\footnote{Such functions  $f$ are called {\em tamed} DTIs in [LapRa\v Zu1--7]; see esp.\ \cite[App.\ A]{fzf} for a development of their general theory.} then $D(f)$, the {\em abscissa of $($absolute $\mathrm{or}$ Lebesgue$)$ convergence} of $f$, is defined as the unique extended real number (i.e., $D(f)\in\eR\cup\{\pm\ty\}$) such that $\{\re s>D(f)\}$ is the {\em maximal} open right half-plane (of the form $\{\re s>\alpha\}$, for some $\alpha\in\eR\cup\{\pm\ty\}$) on which the Lebesgue integral initially defining $f$ in~\eqref{41/4} is convergent (or, equivalently, is absolutely convergent), with $\mu$ replaced by $|\mu|$, the total variation measure of $\mu$.
(Recall that $|\mu|=\mu$ if $\mu$ is positive.)
In short, $D(f)$ is called the {\em abscissa of convergence} of $f$.
Furthermore, the aforementioned maximal right half-plane is denoted by $\Pi(f)$ and is called the {\em half-plane of} (absolute or Lebesgue) {\em convergence} of (the Dirichlet-type integral) $f$.
It is shown in~\cite[{\S}2.1]{fzf} that under mild hypotheses (which are always satisfied in our setting),
$D(f)$ is well defined and (with the notation of~\eqref{41/4} just above) we have, equivalently:\footnote{Let $D:=\ov{\dim}_BA$, for brevity.
In light of Theorem~\ref{an}, for this alternative definition of $D(\zeta_A)$ (or of $D(\widetilde\zeta_A)$), with $A\subseteq\eR^N$ bounded (as in the present situation), it would suffice to restrict oneself to $\alpha\geq 0$ in the right-hand side of~\eqref{41/2}; this follows since $D(\zeta_A)=\ov{\dim}_BA\geq 0$ and (if $D<N$), $D(\zeta_A)=D(\widetilde{\zeta}_A)$. Here, $\widetilde\zeta_A$ stands for the tube zeta function of $A$, defined by Equation~\eqref{tildz}.}
\begin{equation}\label{41/2}
D(f)=\inf\left\{\alpha\in\eR\,:\,\int_E\varphi(x)^{\alpha}\D|\mu|(x)<\ty\right\},
\end{equation}
where (as above) $|\mu|$ is the total variation (local) measure of $\mu$. Under the stated conditions on $f$, we have 
$\Pi(f)\stq\mathcal{H}(f)$; that is, $-\ty\le D_{\rm hol}(f)\le D(f)\le+\ty$.

Note that the distance zeta function $\zeta_A$, defined by \eqref{z}, is a tamed DTI of the form \eqref{41/4}, with $E:=A_\d$, $\f(x):=d(x,A)$ and $\D \mu(x):=d(x,A)^{-N}\D x$. Furthermore, we can clearly take $C:=\d$ in \eqref{tamed}.

\smallskip

The following key result describes some of the basic properties of distance zeta functions.

\begin{theorem}[{[{\rm LapRa\v Zu1,2}]}]\label{an} Let $A$ be an arbitrary bounded subset of $\eR^N$ and let $\delta$ be a fixed positive real number. Then$:$

\medskip

$(a)$ The distance zeta function $\zeta_A$ is holomorphic on $\{\re s>\ov\dim_BA\}$. Moreover, $\Pi(\zeta_A)=\{\re s>\ov\dim_BA\}$; that is,
\begin{equation}
D(\zeta_A)=\ov\dim_BA.
\end{equation}

\medskip

$(b)$ If the box $($or Minkowski$)$ dimension $D:=\dim_BA$ exists, $D<N$  and $\M_*^D(A)>0$, then $\zeta_A(s)\to+\ty$ as $s\in\eR$
converges to $D$ from the right. In particular, $\mathcal{H}(\zeta_A)=\Pi(\zeta_A)=\{\re s>\dim_BA\}$; that is,
\begin{equation}
D_{\rm hol}(\zeta_A)=D(\zeta_A)=\dim_BA.
\end{equation}
\end{theorem}

\medskip

\begin{remark}\label{2.41/2}
$(a)$ It would be of interest to construct (if possible) a class of nontrivial bounded subsets $A$ of $\eR^N$ such that $D_{\rm hol}(\zeta_A)<D(\zeta_A)$. A trivial example is given by $A=[0,1]$, since then $D_{\rm hol}(\zeta_A)=0$ and $D(\zeta_A)=1$.
 
$(b)$ The analog of Theorem \ref{an} holds for the tube zeta function $\widetilde\zeta_A$ (to be introduced in Definition \ref{zeta_tilde} below), except for the fact that in part $(b)$, one no longer needs to assume that $D<N$.
\end{remark}

\medskip

Given a bounded set $A$, it is of interest to know the corresponding poles of the associated distance zeta function $\zeta_A$,  meromorphically extended (if possible) to a neighborhood of the critical line $\{\re s=D(\zeta_A)\}$. Following the terminology of \cite{lapidusfrank12}, these poles are called the {\em complex dimensions} of $A$ and we denote the resulting set of complex dimensions by $\po(\zeta_A)$.\footnote{Strictly speaking, one should talk about the set $\po(\zeta_A,U)$ of {\em visible complex dimensions} relative to a domain $U\stq\Ce$ to which $\zeta_A$ can be meromorphically extended; see [LapRa\v Zu1--4] (along with \cite{lapidusfrank12}). In the examples described in this paper, we have $U:=\Ce$.} We pay particular attention to the set of complex dimensions of $A$ located on the critical line $\{\re s=D(\zeta_A)\}$, which we call  the set of {\em principal complex dimensions} of $A$ and denote by $\dim_{PC} A$.

\medskip

For example, it is well known that for the ternary Cantor set $C^{(1/3)}$, $\dim_BC^{(1/3)}\\ =\log_32$ and, moreover (see \cite[\S1.2.2 and \S2.3.1]{lapidusfrank12}), with $\I:=\sqrt{-1}$,
$$
\dim_{PC}C^{(1/3)}:= \log_32+\frac{2\pi}{\log3}\,\I\Ze.
$$

\medskip

The following result provides an interesting connection between the residue of the distance zeta function of a fractal set at $D:=\dim_BA$ and its Minkowski contents.

\medskip

\begin{theorem}[{[{\rm LapRa\v Zu1,2}]}]\label{pole1}
Assume that $A$ is a bounded subset of $\eR^N$ which is nondegenerate 
$($that is, $0<\M_*^D(A)\le\M^{*D}(A)<\ty$ and, in particular, $\dim_BA=D$$)$, 
and $D<N$. If the distance zeta function $\zeta_A(\,\cdot\,,A_\delta):=\zeta_A$, initially defined by \eqref{z}, can be meromorphically extended\footnote{The existence and construction of meromorphic extensions of fractal zeta functions is discussed in {\S}\ref{merom_ext}. It is studied in a variety of situatons in [LapRa\v Zu1--4,6].} to a neighborhood of $s= D$,
then $D$ is necessarily a simple pole of $\zeta_A(\,\cdot\,,A_\delta)$, and 
\begin{equation}\label{res}
(N-D)\,\M_*^D(A)\le\res(\zeta_A(\,\cdot\,,A_\delta),D)\le(N-D)\,\M^{*D}(A).
\end{equation}
 Furthermore, the value of $\res(\zeta_A(\,\cdot\,,A_\delta), D)$ does not depend on $\delta>0$.
In particular, if $A$ is Minkowski measurable, then 
\begin{equation}\label{pole1minkg1=}
\res(\zeta_A(\,\cdot\,,A_\delta), D)=(N-D)\,\M^D(A).
\end{equation}
\end{theorem}

\medskip

The distance zeta function defined by \eqref{z} is closely related to the tube zeta function of a fractal set, which, in turn, is defined via the tube function $t\mapsto|A_t|$, for $t>0$, of the fractal set $A$, as we now explain.

\medskip

\begin{defn}[{\rm[LapRa\v Zu1,2]}]\label{zeta_tilde}  Let $\delta$ be a fixed positive number, and let $A$ be a bounded subset of $\eR^N$. Then, the {\em tube zeta 
function}
 of $A$, denoted by $\widetilde\zeta_A$, is defined (for all $s\in\Ce$ with $\re s$ sufficiently large) by
\begin{equation}\label{tildz}
\widetilde\zeta_A(s):=\int_0^\delta t^{s-N-1}|A_t|\,\D t.
\end{equation}
\end{defn} 

\medskip

For any fixed positive real number $\d>0$, the distance and tube zeta functions associated with a given fractal set $A$ are connected as follows:\footnote{We write here $\zeta_A(\,\cdot\,,A_{\delta}):=\zeta_A$ and $\widetilde\zeta_A(\,\cdot\,,\d):=\widetilde\zeta_A$, for emphasis.}
\begin{equation}\label{equ_tilde}
\zeta_A(s,A_\delta)=\delta^{s-N}|A_\delta|+(N-s)\widetilde\zeta_A(s,\d),
\end{equation}
for $\re s>\ov\dim_BA$;\footnote{In light of the principle of analytic continuation, one deduces that identity~\eqref{equ_tilde} continues to hold whenever one (and hence, both) of the fractal zeta functions $\zeta_A$ and $\widetilde\zeta_A$ is meromorphic on a given domain $U\subseteq\Ce$. }
see [LapRa\v Zu1,2].\footnote{The case when $D=N$ in Theorem \ref{pole1mink_tilde} must be treated separately.}
Using this result, it is easy to 
obtain the analog of Theorem \ref{an} for $\widetilde\zeta_A$ (as was stated in Remark \ref{2.41/2}$(b)$ above) and to
reformulate Theorem \ref{pole1} in terms of the tube zeta functions.
In particular, we conclude that the residue of the tube zeta function of a fractal set, computed at $s=D$, is equal to
 its Minkowski content, provided the set is Minkowski measurable.

\medskip

\begin{theorem}[{[{\rm LapRa\v Zu1,2}]}]\label{pole1mink_tilde}
Assume that $A$ is a nondegenerate bounded subset of $\eR^N$ $($so that $D:=\dim_BA$ exists$)$, 
and there exists a meromorphic extension of $\widetilde\zeta_A$ to a neighborhood of $D$. Then,  $D$ is a simple  pole of $\widetilde\zeta_A$,
and for any positive $\delta$, $\res(\widetilde{\zeta}_A,D)$ is independent of $\delta$. Furthermore, we have
\begin{equation}\label{zeta_tilde_M}
\M_*^D(A)\le\res(\widetilde\zeta_A, D)\le \M^{*D}(A).
\end{equation}
In particular, if $A$ is Minkowski measurable, then 
\begin{equation}\label{zeta_tilde_Mm}
\res(\widetilde\zeta_A, D)=\M^D(A).
\end{equation}
\end{theorem}

\medskip

A class of fractal sets $A$ for which we have strict inequalities in (\ref{zeta_tilde_M}) (and hence also in~\eqref{res} of Theorem~\ref{pole1} above) is constructed in
Theorem \ref{nonmeasurable}; see (\ref{res_inequalities}).

\medskip

\section{Meromorphic extensions of fractal zeta functions}\label{merom_ext}

Since the definition of the set of principal complex dimensions $\dim_{PC} A$ of $A$ requires  the existence of a suitable meromorphic extension of the distance zeta function $\zeta_A$,
 it is natural to study this issue in more detail.
For simplicity, we formulate the results of this section for $\widetilde\zeta_A$, but we note that the analogs of Theorems \ref{measurable} and \ref{nonmeasurable} also hold for $\zeta_A$,  provided $D<N$; see \cite[{\S}2.3.3]{fzf} or \cite{mezf}.

\medskip

\begin{theorem}[Minkowski measurable case, {\rm[LapRa\v Zu1,3]}]\label{measurable}%
Let $A$ be a bounded subset of $\eR^N$ such that there exist  $\alpha>0$, $\mathcal {M}\in(0,+\ty)$ and $D\ge0$ satisfying
\begin{equation}\label{A_t}
|A_t|= t^{N-D}\left({\mathcal {M}}+O(t^\alpha)\right)\quad\mathrm{as}\quad t\to0^+.
\end{equation}
Then, $\dim_BA$ exists and $\dim_BA=D$. Furthermore, $A$ is Minkowski measurable with Minkowski content $\mathcal {M}^D(A)=\mathcal {M}$. 
Moreover, the tube zeta function $\widetilde\zeta_A$ has for abscissa of convergence $D(\widetilde\zeta_A)=\dim_BA=D$ and possesses a $($necessarily unique$)$ meromorphic continuation $($still denoted by $\widetilde\zeta_A)$ to $($at least$)$ the open right half-plane 
$\{\re s>D-\alpha\}$.
The only pole of $\widetilde\zeta_A$ in this half-plane is $s=D$; it is simple and, moreover, $\res(\widetilde\zeta_A,D)=\M$. 
\end{theorem}

\medskip

Next, we deal with a useful class of Minkowski nonmeasurable sets. 
Before stating Theorem \ref{nonmeasurable}, let us first introduce some notation.
Given a locally integrable $T$-periodic function $G:\eR\to\eR$, with $T>0$, we denote by $G_0$ its truncation to $[0,T]$, while the Fourier transform of $G_0$
is denoted by $\hat G_0$: for all $t\in\eR$,
\begin{equation}\label{fourier}\noindent
\hat G_0(t):=\int_{-\ty}^{+\ty}{\E}^{-2\pi {\I}t\tau}G_0(\tau)\,\D\tau=\int_0^T{\E}^{-2\pi {\I}t\tau}G(\tau)\,\D\tau.
\end{equation}

\medskip

\begin{theorem}[Minkowski nonmeasurable case,  {\rm[LapRa\v Zu1,3]}]\label{nonmeasurable}%
Let $A$ be a bounded subset of $\eR^N$ such that there exist $D\ge0$, $\alpha>0$, and $G:\eR\to(0,+\ty)$ a nonconstant periodic function with period $T>0$, 
satisfying
\begin{equation}\label{G}
|A_t|=t^{N-D}\left(G(\log t^{-1})+O(t^\alpha)\right)\quad\mbox{{\rm as}\qs$t\to0^+$.}
\end{equation}
 Then $G$ is continuous, $\dim_BA$ exists and $\dim_BA=D$. Furthermore, $A$ is Minkowski nondegenerate,
 with upper and lower Minkowski contents respectively given by
\begin{equation}\label{1.4.201/2}
\M_*^D(A)=\min G,\quad \M^{*D}(A)=\max G.
\end{equation}
Moreover, the tube zeta function $\widetilde\zeta_A$ has for abscissa of convergence $D(\widetilde\zeta_A)=D$ and possesses a $($necessarily unique$)$ meromorphic extension $($still denoted by $\widetilde\zeta_A$$)$
to {\rm({\it at least})} the half-plane $\{\re s>D-\alpha\}$.

\medskip

In addition, the set of principal complex dimensions of $A$ is given by
\begin{equation}\label{Dpoles}
\dim_{PC} A=\left\{s_k=D+\frac{2\pi}T{\I}k:\hat G_0\Big(\frac kT\Big)\ne0,\,\,k\in\Ze\right\}
\end{equation}
{\rm({\it see} (\ref{fourier}))} and there are no other complex dimensions in $\{\re s>D-\alpha\}$; they are all simple, and the residue at each $s_k\in\mathcal \dim_{PC} A$, with $k\in\Ze$, is given by
\begin{equation}\label{res_fourier}
\res(\widetilde\zeta_A,s_k)=\frac1T\hat G_0\Big(\frac kT\Big).
\end{equation}
If $s_k\in \dim_{PC} A$, then $s_{-k}\in \dim_{PC} A$ $($in agreement with the `reality principle'$)$, and
$|\res(\widetilde\zeta_A,s_k)|\le \frac1T\int_0^TG(\tau)\,\D\tau$; 
furthermore, $\lim_{k\to\pm\ty}\res(\widetilde\zeta_A,s_k)=0$.

\medskip

Moreover, the set of principal complex dimensions of $A$ contains $s_0=D$, and
\begin{equation}\label{avarage}
\res(\widetilde\zeta_A,D)=\frac1T\int_0^TG(\tau)\,\D\tau.
\end{equation}
In particular, $A$ is {\rm not} Minkowski measurable and
\begin{equation}\label{res_inequalities}
\M_*^D(A)<\res(\widetilde\zeta_A,D)<\M^{*D}(A).
\end{equation}
\end{theorem}

\medskip

\begin{example}[{$a$-{strings}}]\label{a-string2}
The compact set $A:=\{j^{-a}:j\in\eN\}\cup\{0\}$, where $a>0$, is Minkowski measurable and
\begin{equation}\label{a-string}
\M^D(A)=\frac{2^{1-D}}{1-D}a^D,\quad D:=\dim_BA=\frac1{1+a}.
\end{equation} 
(See $[$Lap 1, Exple.\ 5.1 and App.\ C$]$.)
The associated fractal string $\mathcal {L}=(\ell_j)_{j\ge1}$, defined by $\ell_j=j^{-a}-(j+1)^{-a}$ for all $j\ge1$ (or, equivalently, by $\O:=[0,1]\setminus A\st\eR$, so that $\pa\O=A$), is called the {\em $a$-string}; see [Lap1--3], [LapPo1--2], [HeLap] and \cite[{\S}6.5.1]{lapidusfrank12}. In light of (\ref{pole1minkg1=}) and (\ref{zeta_tilde_Mm}), we then know that
$\res(\zeta_A(\,\cdot\,,A_\delta),D)=(1-D)\M^D(A)$ and $\res(\widetilde\zeta_A,D)=\M^D(A)$.
\end{example}

\begin{example}[{{fractal nests}}]
Let $a>0$ and let $A$ be the countable union of concentric circles in $\eR^2$, centered at the origin and of radii $r=k^{-a}$,  where $k\in\eN$. According to the terminology introduced in [LapRa\v Zu1--4], $A$ is called the {\em fractal nest of inner type} generated by the $a$-string from the preceding example.  Then, using the the distance zeta function of $A$ it is possible to show that 
\begin{equation}\label{dim_nest}
D:=\ov\dim_BA=\max\Big\{1,\frac2{1+a}\Big\}.
\end{equation}
(See \cite[Chap.\ 3]{fzf} and [LapRa\v Zu2--4].) The set $A$ is closely related to the planar spiral $\C$ defined in polar coordinates by $r=\theta^{-a}$, $\theta\ge \theta_0$, where $\theta_0>0$, and the value of $\dim_B\C$ is the same as for $A$; see \cite{tricot}. We mention in passing that for $a\ne1$, the fractal nest $A$ (as well as the corresponding spiral $\C$) is Minkowski measurable and for every $a\in(0,1)$, the value of its Minkowski content is independent of $\theta_0$ and given by
\begin{equation}\label{AGamma}
\M^D(A)=\pi(2/a)^{2a/(1+a)}\frac{1+a}{1-a}.
\end{equation}
Using \eqref{AGamma}, along with with Eq.\  \eqref{pole1minkg1=} from Theorem \ref{pole1}, we conclude that the residue of the distance zeta function $\zeta_A$, as well as of $\zeta_\C$, computed at $s=D$, is given by
\begin{equation}
\res(\zeta_A,D)=\res(\zeta_\C,D)=\pi(2/a)^{2a/(1+a)}\frac{2a}{1-a},
\end{equation}
provided $a\in(0,1)$. For $a=1$, we have $\M^1(A)=\M^1(\C)=+\ty$. These and related results are useful in the study of fractal properties of spiral trajectories of planar vector fields; see, e.g., \cite{zuzu}. 

More generally, if we consider the fractal nest $A_N$ defined as the countable union of concentric spheres in $\eR^N$, centered at the origin and of radii $r=k^{-\a}$, where $k\in\eN$, then using the distance zeta function $\zeta_{A_N}$, it can be shown (see \cite[{\S}3.4]{fzf} and [LapRa\v Zu2--4]) that
\begin{equation}
\ov\dim_B A_N=\max\Big\{N-1,\frac N{1+a}\Big\}. 
\end{equation}
Note that for $N=1$ and $N=2$, we recover the box dimension of the $a$-string and of the fractal nest, respectively; see  [LapRa\v Zu2--4] and  Eqs.\ \eqref{a-string}--\eqref{dim_nest} above.
\end{example}

\medskip

In the following result, we provide the distance zeta function of the Sierpi\'nski carpet and the corresponding principal complex dimensions. It is well known that the Sierpi\'nski carpet is not Minkowski measurable. See, e.g., \cite{lapidusfrank12}, as well as \cite{lana} for explicit values of its upper and lower Minkowski contents. A similar result can be obtained for the Sierpi\'nski gasket (and its higher-dimensional analogs); see \cite[{\S}3.2.2]{fzf} and [LapRa\v Zu2--4].

\begin{prop}[Distance zeta function of the Sierpi\'nski carpet]\label{sierpinski_carpet0}
Let $A$ be the Sierpi\'nski carpet in $\eR^2$, constructed in the usual way inside the unit square. Let $\d$ be a fixed positive real number. 
We assume without loss of generality that $\d>1/6$ $($so that for this choice of $\d$, $A_\d$ coincides with the $\d$-neighborhood of the unit square $[0,1]^2$$)$.
Then, for all $s\in\Ce$, the distance zeta function $\zeta_A$ of the Sierpi\'nski carpet is given by
\begin{equation}\label{zeta_carpet}
\zeta_A(s)=\frac{8}{2^ss(s-1)(3^s-8)}+2\pi\frac{\d^s}s+4\frac{\d^{s-1}}{s-1},
\end{equation}
which is meromorphic on the whole complex plane. In particular, the set of complex dimensions and of principal complex dimensions of the Sierpi\'nski carpet are given, respectively, by
\begin{equation}
\po(\zeta_A)=\{0,1\}\cup\dim_{PC} A,\q\dim_{PC} A= \log_38+\frac{2\pi}{\log 3}\I\Ze.
\end{equation}
Furthermore, each of the complex dimensions $($i.e., each of the poles of $\zeta_A$$)$ is simple. Moreover, the residues of the distance zeta function $\zeta_A$ computed
at the principal poles $s_k:=\log_38+\frac{2\pi}{\log 3}k\I$, with $k\in\Ze$, are given by
\begin{equation}\label{ressk}
\res(\zeta_A,s_k)=\frac{2^{-s_k}}{(\log3)s_k(s_k-1)}.
\end{equation}
Finally, the approximate values of the lower and upper $D$-dimensional Minkowski contents are given by 
$\M_*^D(A)\approx 1.350670$ and $\M^{*D}(A)\approx1.355617$.
$($The precise values can be found in \cite{lana}$.)$
\end{prop}

\smallskip

\begin{proof}[Sketch of the proof]
In order to compute the distance zeta function 
$$
\zeta_A(s):=\int_{A_\d}d((x,y),A)^{s-2}\D x\,\D y
$$
of the Sierpi\'nski carpet $A$, we first have to calculate 
\begin{equation}\label{rfd}
\zeta_{A_k}(s,\O_k):=\int_{\O_k}d((x,y),A_k)\,\D x\,\D y, 
\end{equation}
where $\O_k$ is a square of the $k$-th generation (of side lengths $a_k=3^{-k}$) and $A_k$ is its boundary. (Here, we deal in fact with `relative distance zeta functions', which are discussed in Remark \ref{relative} just below; see \eqref{relativezf} and [LapRa\v Zu1--4].) 
This can be easily done by splitting $\O_k$ into the disjoint union of eight congruent right-angle triangles, and we obtain after a short computation that
$
\zeta_{A_k}(s,\O_k)=8\cdot 2^{-s}a_k^ss^{-1}(s-1)^{-1}
$.
Since the $k$-th generation consists of $8^{k-1}$ squares congruent to $\O_k$, we deduce that
\begin{equation}\label{Akv}
\zeta_A(s,[0,1]^2)=\sum_{k=1}^\ty8^{k-1}\zeta_{A_k}(s,\O_k)=\frac{8}{2^ss(s-1)(3^s-8)},
\end{equation}
for $\re s>\log_38$. 
The last expression in \eqref{Akv} is meromorphic in all of $\Ce$. Hence, upon analytic continuation, $\zeta_A(s,[0,1]^2)$ is given by that expression for all $s\in\Ce$. Note that
the value of $\zeta_A(s,[0,1]^2)$ is precisely equal to the first term on the right-hand side of \eqref{zeta_carpet}. The remaining two terms are obtained by considering $\zeta_A(s,A_\d\setminus[0,1]^2)$, which can be easily reduced to considering a disk $B_{\d}(0)$ of radius $\d$ with respect to its origin $0\in\eR^2$, and two rectangles that are congruent to $\O_0:=(0,1)\times(-\d,\d)$ with respect to its middle section $A_0:=(0,1)\times\{0\}$.
\end{proof}

\begin{remark}\label{relative}
Eq.\  \eqref{rfd} is a very special case of the zeta function of a {\em relative fractal drum} $(A,\O)$ in $\eR^N$, a notion which will be briefly discussed in {\S}\ref{fzfrfd} and is the object of [LapRa\v Zu4] and \cite[Chap.\ 4]{fzf}; see the first equality in Eq.\  \eqref{relativezf} below.
\end{remark}
\medskip

\section{Transcendentally quasiperiodic sets}\label{qp_sets}

In this section, we define a class of {\em quasiperiodic fractal sets}. The simplest of such sets has two incommensurable periods. Moreover, using suitable generalized Cantor sets, it is possible to ensure that
the quotient of their periods be a transcendental real number.
Our construction of such sets is based on a class of generalized Cantor sets with two parameters, which we now introduce.

\medskip

\begin{defn}[{\rm[LapRa\v Zu1,2]}]\label{Cma}
The generalized Cantor sets $C^{(m,a)}$ are determined by an integer $m\ge2$ and a real number $a\in(0,1/m)$.
In the first step of the analog of Cantor's construction, we start with $m$ equidistant, closed intervals in $[0,1]$ of length $a$, with $m-1$ `holes', each of length $(1-ma)/(m-1)$. In the second step, we continue by scaling by the factor $a$ each of the $m$ intervals of length $a$; and so on, ad infinitum.
The  $($two-parameter$)$ {\em generalized Cantor set} $C^{(m,a)}$ is then defined as the intersection of the decreasing sequence of compact sets constructed in this way.
It is easy to check that $C^{(m,a)}$ is a perfect, uncountable compact subset of $\eR$; furthermore, $C^{(m,a)}$ is also self-similar.
For $m=2$, the sets $C^{(m,a)}$ are denoted by  $C^{(a)}$.
The classic ternary Cantor set is obtained as $C^{(2,1/3)}$.
In order to avoid any possible confusion, we note that the generalized Cantor sets introduced here are different from the generalized Cantor strings introduced and studied in~\cite[Chap.\ 10]{lapidusfrank12}, as well as used in a key manner in \cite[Chap.\ 11]{lapidusfrank12}.
\end{defn}

\smallskip

We collect some of the basic properties of generalized Cantor sets in the following proposition.

\smallskip

\begin{prop}[Generalized Cantor sets, {\rm[LapRa\v Zu1,2]}]\label{Cmap}
 If $A:=C^{(m,a)}$ is the generalized Cantor set introduced in Definition~\ref{Cma}, 
where $m$ is an integer larger than $1$, and $a\in(0,1/m)$,
then
\begin{equation}\label{2.1.1}
D:=\dim_BA=D(\zeta_A)=\log_{1/a}m.
\end{equation}
Furthermore, the tube formula associated with $A$ is given by
\begin{equation}\label{Cmat}
|A_t|=t^{1-D}G(\log t^{-1})\q\mbox{\rm for all\q $t\in(0,t_0)$,}
\end{equation}
where $t_0$ is a suitable positive constant and $G=G(\tau)$ is a continuous, positive and nonconstant periodic function, with minimal period $T=\log (1/a)$.

\medskip

Moreover, $A$ is Minkowski nondegenerate and Minkowski nonmeasurable; that is, $0<\mathcal{M}_*^D(A)<\mathcal{M}^{*D}(A)<\ty$.\footnote{The periodic function $G=G(\tau)$, as well as the values of $\mathcal{M}_*^D(A)$ and $\mathcal{M}^{*D}(A)$, can be explicitly computed; see [LapRa\v Zu1, {\S}3.1.1]}

\smallskip

Finally, the distance zeta function of $A$ admits a meromorphic continuation to all of $\Ce$ and
the set of principal complex dimensions of $A$ is given by
\begin{equation}\label{2.1.6}
\dim_{PC}A=D+\frac{2\pi}T{\I}\Ze.
\end{equation}
Besides $(\dim_{PC}A)\cup\{0\}$, there are no other poles, and all of the poles of $\zeta_A$ are simple.
In particular, $\po(\zeta_A)=(D+\frac{2\pi}T{\I}\Ze)\cup\{0\}$.
\end{prop}

\smallskip

The definition of quasiperiodic sets is based on the following notion of quasiperiodic functions, which will be useful for our purposes.\footnote{We note that Definition~\ref{quasip}, although rather close to the one provided in \cite{enc}, is very different from the usual definition of Bohr-type quasiperiodic functions.}

\smallskip

\begin{defn}[{\rm[LapRa\v Zu1,2]}]\label{quasip}
We say that a function $G=G(\tau):\eR\to\eR$ is {\em transcen\-dentally 
$n$-quasiperiodic}, with $n\ge2$,
 if it is of the form $G(\tau)=H(\tau,\dots,\tau)$,
where $H:\eR^n\to\eR$  is a function that is nonconstant and $T_k$-periodic in its $k$-th component, for each $k=1,\dots,n$, and the periods $T_1,\dots, T_n$ are algebraically (and hence, rationally) independent.\footnote{That is, linearly independent over the field of algebraic numbers.} The positive numbers $T_i$ ($i=1,\dots,n$) are called the 
{\em quasiperiods} of $G$. If, instead, the set of quasiperiods $\{T_1,\dots,T_n\}$ is rationally independent and algebraically dependent, we say that {\em $G$ is algebraically $n$-quasiperiodic}.
\end{defn}

\begin{defn}[{[LapRa\v Zu1,2]}]\label{quasiperiodic}
Given a bounded subset $A$ of $\eR^N$, we say that a function $G:\eR\to\eR$ {\em is associated with the set A} (or {\em corresponds to $A$}) if it is nonnegative and $A$ has the following tube formula:
\begin{equation}
|A_t|=t^{N-D}(G(\log(1/t))+o(1))\textrm{ as }t\to0^+,
\end{equation}
with $0<\liminf_{\tau\to+\ty} G(\tau)\le\limsup_{\tau\to+\ty} G(\tau)<\ty$.
In addition, we say that $A$ is a {\em transcendentally $(${\rm resp.}, algebraically$)$ $n$-quasiperiodic set} if the function $G=G(\tau)$ is transcendentally $($resp., algebraically$)$
$n$-quasiperiodic. The smallest possible value of $n$ is called the {\em order of quasiperiodicity} of $A$.
\end{defn}

The following result, which has a variety of generalizations as will be briefly explained below, provides a construction of transcendentally $2$-quasi\-peri\-odic fractal sets. Its proof is based on the classical Gel'fond--Schneider theorem (as described in~\cite{gelfond}) from transcendental number theory.

\begin{theorem}[{[{\rm LapRa\v Zu1,2}]}]\label{trans} Let $C^{(m_1,a_1)}$ and $\dim_B C^{(m_2,a_2)}$ be two generalized Cantor sets such that their box dimensions coincide and are equal to $D\in(0,1)$. Assume that $I_1$ and $I_2$ are two unit closed intervals of $\eR$, with disjoint interiors,
and define $A_1:=(\min I_1)+C^{(m_1,a_1)}\st I_1$ and $A_2:=(\min I_2)+C^{(m_2,a_2)}\st I_2$. 
Let $\{p_1,p_2,\dots,p_k\}$ be the set of all distinct prime factors of $m_1$ and $m_2$, and write
\begin{equation}\nonumber
m_1=p_1^{\alpha_1}p_2^{\alpha_2}\dots p_k^{\alpha_k},\quad m_2=p_1^{\beta_1}p_2^{\beta_2}\dots p_k^{\beta_k},
\end{equation}
where $\alpha_i,\beta_i\in\eN\cup\{0\}$\label{n_0} for $i=1,\ldots,k$. If the exponent vectors $\mathbf{e}_1$ and $\mathbf{e}_2$ of, respectively, $m_1$ and $m_2$, defined by
\begin{equation}\nonumber
\mathbf{e}_1:=(\alpha_1,\alpha_2,\dots,\alpha_k)\,\,\qs\mathrm{and}\qs\,\,\mathbf{e}_2:=(\beta_1,\beta_2,\dots,\beta_k),
\end{equation}
are linearly independent over the field of rational numbers, then the compact set $A:=A_1\cup A_2\st\eR$ is transcendentally $2$-quasiperiodic.

\smallskip

Moreover, the distance zeta function $\zeta_A$ can be meromorphically extended to the whole complex plane, and we have that
$D(\zeta_{A})=D$. The set $\dim_{PC}A$ of principal complex dimensions of $A$ is given by
\begin{equation}\label{511/2}
\dim_{PC} A=
D+\Big(\frac{2\pi}{T_1}\Ze\cup \frac{2\pi}{T_2}\Ze \Big)\I.
\end{equation}
Besides $(\dim_{PC}A)\cup\{0\}$, there are no other poles of the distance zeta function~$\zeta_{A}$ and they are all simple.
In particular, 
\begin{equation}
\po(\zeta_A)=\Bigg(D+\Big(\frac{2\pi}{T_1}\Ze\cup \frac{2\pi}{T_2}\Ze \Big)\I \Bigg)\cup\{0\}.
\end{equation}
\end{theorem}

\begin{remark}\label{tyqp} This result can be considerably extended by using Baker's theorem \cite[Thm.\ 2.1]{baker} which, in turn, is a
far-reaching extension of the aforementioned Gel'fond--Schneider's theorem.
Indeed, for any fixed integer $n\ge2$, using Baker's theorem and $n$ generalized Cantor sets,
an explicit construction of a class of transcendentally $n$-quasiperiodic fractal sets is  given in [LapRa\v Zu2] and \cite[\S3.1]{fzf}. 
In [LapRa\v Zu4] and \cite[Chap.\ 4]{fzf}, we even construct a set which is transcendentally $\infty$-quasiperiodic; see \S\ref{maxhf}.
\end{remark}

\section{Maximally hyperfractal $\ty$-quasiperiodic sets}\label{maxhf}

It is possible to construct a bounded subset $A$ of the real line, such that the corresponding distance zeta function $\zeta_A$ has for abscissa of (Lebesgue, i.e., absolute) convergence $D(\zeta_A)$ any prescribed real number $D\in(0,1)$ and $A$ is {\em maximally hyperfractal}; that is, any point on the critical line $\{\re s=D\}$ is a nonremovable singularity of the corresponding distance zeta function $\zeta_A$. In particular, there is no meromorphic continuation of $\zeta_A$ to any open and connected neighborhood of the critical line (and, moreover, not even to any open and connected neighborhood of an arbitray point on the critical line). Furthemore, it is possible to construct a maximally hyperfractal set which is $\ty$-transcendentally quasiperiodic as well. A construction of such sets is described in detail in \cite[Chap.\ 4]{fzf} or in \cite{rfds}.
In the sequel, we provide a rough sketch of this construction.

The set $A\st\eR$ which is a maximal hyperfractal and $\ty$-transcendentally quasiperiodic set, can be constructed as the nonincreasing sequence 
\begin{equation}
A=A_{\mathcal L}=\bigg\{a_k:=\sum_{j=k}^\ty \ell_j:k\in\eN\bigg\}
\end{equation}
of positive real numbers $a_k$ converging to zero as $k\to\ty$, generated by a suitable bounded fractal string $\mathcal L=(\ell_j)_{j\ge1}$. Roughly speaking, the fractal string $\mathcal L$ is obtained as a (suitably defined) {\em union of an infinite sequence of bounded fractal strings} ${\mathcal L}_k:=(\ell_{kj})_{j\ge1}$, corresponding to generalized Cantor sets of the form $c_k\cdot C^{(m_k,a_k)}$, for $k\in\eN$, with carefully chosen values of the parameters $m_k$ and $a_k$ appearing in Definition \ref{Cma}, and where $(c_k)_{k\ge1}$ is an appropriate summable sequence of positive real numbers. 

More precisely, the union 
${\mathcal L}:=\bigsqcup_{k=1}^\ty {\mathcal L}_k$ 
of the sequence of bounded fractal strings ${\mathcal L}_k$ is defined as the set-theoretic union of the elements of the strings, but by definition, each of its elements has for multiplicity the sum of the corresponding multiplicities from all of the fractal strings ${\mathcal L}_k$ to which belongs the element in question. Note that the multiplicity of an element of $\mathcal L$ is well defined since this element must belong to at most finitely many bounded fractal strings ${\mathcal L}_k$, which follows from the fact that the sequence $c_k$ converges to $0$ as $k\to\ty$. Moreover, we must assume that $\sum_{k=1}^\ty c_k<\ty$, so that the string $\mathcal L$ be bounded (i.e., $\sum_{j=1}^\ty\ell_j<\ty$). We can also ensure that for each positive integer $k$, the corresponding upper box dimension of ${\mathcal L}_k$ (that is, of the set $A_{{\mathcal L}_k}$) be equal to a fixed value of $D\in(0,1)$, prescribed in advance. (Note that the set $A_{\mathcal L}$ is distinct from $\cup_{k=1}^\ty A_{{\mathcal L}_k}$.)

Recall that the oscillatory period of ${\mathcal L}_k$ (in the sense of \cite{lapidusfrank12}), which is defined by
${\mathbf p}_k:=\frac{2\pi}{\log(1/a_k)}$,
provides valuable information about the density of the set of principal complex dimensions of ${\mathcal L}_k$ on the critical line $\{\re s=D\}$. More precisely,
by choosing the coefficient $a_k\in(0,1/m_k)$ so that $a_k\to0$ as $k\to\ty$, we see that for the set of principal complex dimensions of the generalized Cantor string ${\mathcal L}_k$ (i.e, the set of the principal poles of $\zeta_{{\mathcal L}_k}$),
$\dim_{PC}{\mathcal L}_k=\dim_{PC} C^{(m_k,a_k)}=D+{\mathbf p}_k\I\Ze$,
becomes denser and denser on the critical line, as $k\to\ty$, since then  the oscillatory period ${\mathbf p}_k$ tends to zero. Therefore, the distance zeta function of the fractal string ${\mathcal L}:=\sqcup_{k=1}^\ty{\mathcal L}_k$ will have 
$D+\Big(\bigcup_{k=1}^\ty{\mathbf p}_k\Ze\Big)\I$
as a set of singularities, which is densely packed on the critical line 
$\{\re s=D\}=D+\eR\I$,  
since the set $\cup_{k=1}^\ty{\mathbf p}_k\Ze$ is clearly dense in $\eR$.
In conclusion, the whole critical line $\{\re s=D\}$ consists of nonremovable singularities of $\zeta_{\mathcal L}$,\footnote{In light of the discussion surrounding Eq.\  \eqref{dist_geo} above, the same is true if $\zeta_{\mathcal L}$ is replaced by $\zeta_{A_{\mathcal L}}$ or, more generally, by the {\em relative distance zeta function} $\zeta_A(\,\cdot\,,\O)$ defined by $\zeta_A(s,\O):=\int_\O d(x,A)^{s-N}\D x$ (see {\S}\ref{fzfrfd} and \cite{rfds} or \cite[Chap.\ 4]{fzf}), where $A=\pa\O$ is the boundary of any geometric realization of $\mathcal L$ by a bounded open subset $\O$ of $\eR$.} which by definition means that the fractal string $\mathcal L$ is maximally hyperfractal. Hence, the corresponding set $A:=A_{\mathcal L}$ is also maximally hyperfractal.

Since the coefficients $a_k$, appearing in the definition of the generalized Cantor set (see Definition \ref{Cma}), have been chosen above so that $a_k\to0$ as $k\to\ty$, it is clear that  $m_k\to\ty$, because
$D=\dim_B C^{(m_k,a_k)}=\frac{m_k}{\log(1/a_k)}$,
where $D\in(0,1)$ is given in advance and independent of $k$. This enables us to use our result mentioned in Remark \ref{tyqp}, obtained by means of Baker's theorem from transcendental number theory [Ba],
in order to ensure that the sequence of quasiperiods 
$T_k:=\log(1/a_k)$, $k\in\eN$,
is algebraically independent (that is, any finite subset of this set of quasiperiods is linearly independent over the field of algebraic real numbers).\footnote{The algebraic independence of the set of quasiperiods $\{T_k:k\ge1\}$, with $k\ge1$, can be deduced (using the aforementioned Baker's theorem, [Ba]) if we assume, in addition, that the sequence $(\mathbf{e}_k)_{k\ge1}$ (suitably redefined), corresponding to the sequence $(m_k)_{k\ge1}$, is rationally independent.} According to Definition  \ref{quasip}, this means that $\mathcal L$ is $\ty$-transcendentally quasiperiodic, and so is the corresponding bounded subset $A:=A_{\mathcal L}$ of the real line.

As we see from the above rough description, the nature of a subset $A:=A_{\mathcal L}$ of the real line which is maximally hyperfractal and $\ty$-transcendentally quasiperiodic, is in general extremely complex, although it is, in fact, `just' defined in terms of a nonincreasing sequence of positive real numbers converging to zero.

In closing this discussion, we mention that this construction (as well as Theorem \ref{trans} and its generalization mentioned in Remark \ref{tyqp}), extends to any $N\ge2$, by letting $B:=A\times[0,1]^{N-1}\st\eR^N$; see [LapRa\v Zu1,2,4].

\section{Fractal zeta functions of relative fractal drums}\label{fzfrfd}

In this section, we survey some of the definitions and results from \cite[Chap.\ 4]{fzf}; see also \cite{rfds}.
Let $A$ be a (possibly unbounded) subset of $\eR^N$ and let $\O$ be a (possibly unbounded) Borel subset of $\eR^N$ of finite $N$-dimensional Lebesgue measure. We say that the ordered pair $(A,\O)$ is a {\em relative fractal drum} (or RFD, in short) if there exists a positive real number $\d$ such that $\O\stq A_\d$. It is easy to see that for every $\d>0$, any bounded subset $A$ can be identified with the relative fractal drum $(A,A_\d)$. Furthermore, any bounded fractal string $\mathcal L=(\ell_j)_{j=1}^\ty$ can be identified with the relative fractal drum $(\cup_{j=1}^\ty\pa I_j,\cup_{j=1}^\ty I_j)$, where $(I_j)_{j=1}^\ty$ is a family of pairwise disjoint open intervals in $\eR$ such that $|I_j|_1=\ell_j$ for all $j\ge1$.

Given a relative fractal drum $(A,\O)$ in $\eR^N$ and for a fixed real number $r$, we define  the {\em relative upper} and {\em relative lower $r$-dimensional Minkowski contents} of $(A,\O)$, respectively, by
$$
\M^{*r}(A,\O):=\limsup_{t\to0^+}\frac{|A_t\cap\O|}{t^{N-r}},\q \M_*^r(A):=\liminf_{t\to0^+}\frac{|A_t\cap\O|}{t^{N-r}}.
$$
The {\em relative upper} and {\em relative lower box $(\mathrm{or}$ Minkowski$)$ dimensions} of $(A,\O)$ are then defined, respectively, by
\begin{equation}\label{rbdim}
\begin{aligned}
\ov\dim_B(A,\O)&:=\inf\{r\in\eR:\M^{*r}(A,\O)=0\},\\
\underline\dim_B(A,\O)&:=\inf\{r\in\eR:\M_*^r(A,\O)=0\}.
\end{aligned}
\end{equation}
It is easy to check that $-\ty\leq\underline{\dim}_B(A,\O)\leq\ov{\dim}_B(A,\O)\leq N$, and it is shown in [LapRa\v Zu1,4] that the relative box dimensions can indeed attain arbitrary negative values as well, including $-\ty$ (an obvious example is when $A_\d\cap\O=\emptyset$ for some $\d>0$). Intuitively, negative relative box dimensions correspond to the {\em property of flatness} of the RFD under consideration.  If $\dim_B(A,\O)=-\ty$, then the RFD $(A,\O)$ is said to be {\em infinitely flat}. A nontrivial example of an infinitely flat RFD $(A,\O)$ in $\eR^2$ is given by $A:=\{(0,0)\}$ and $\O:=\{(x,y)\in(0,1)^2:0<y<\E^{-1/x}\}$. Other examples of flat RFDs can be found in~[LapRa\v Zu1,4].

If $(A,\O)$ is such that $\underline\dim_B(A,\O)=\ov\dim_B(A,\O)$, then this common value is denoted by $\dim_B(A,\O)$ and is called the {\em box $(\mathrm{or}$ Minkowski$)$ dimension of $(A,\O)$}.
Moreover, if $(A,\O)$ is such that, for some $D\in(-\ty,N]$, we have $0<\M_*^D(A,\O)\le\M^{*D}(A,\O)<\ty$ (in particular, then $\dim_B(A,\O)$ exists and $D=\dim_B(A,\O)$), we say that {\em $(A,\O)$ is  Minkowski nondegenerate}. If $\M_*^D(A,\O)=\M^{*D}(A,\O)$, then the common value is denoted by $\M^D(A,\O)$ and called the {\em Minkowski content of} $(A,\O)$. Finally, assuming that $(A,\O)$ is such that $\M^D(A,\O)$ exists and $0<\M^D(A,\O)<\ty$, we say that the RFD {\em $(A,\O)$ is Minkowski measurable}.

To any given RFD $(A,\O)$ in $\eR^N$, we can associate the corresponding {\em relative distance zeta function} and the {\em relative tube zeta function} defined, respectively, by
\begin{equation}\label{relativezf}
\zeta_A(s,\O):=\int_\O d(x,A)^{s-N}\D x,\q
\widetilde\zeta_A(s,\O):=\int_0^\d t^{s-N-1}|A_t\cap\O|\,\D t,
\end{equation}
for all $s\in\Ce$ with $\re s$ sufficiently large, where $\d$ is a fixed positive real number.
They are a valuable theoretical and technical new tool in the study of fractals. 

The basic result dealing with relative distance zeta functions, analogous to Theorem \ref{an} of {\S}\ref{disttube}, is provided by the following theorem.

\begin{theorem}[{[{\rm LapRa\v Zu1,4}]}]\label{an2}
 Let $(A,\O)$ be an arbitrary RFD. Then$:$

\medskip

$(a)$ The distance zeta function $\zeta_A(\,\cdot\,,\O)$ is holomorphic on $\{\re s>\ov\dim_B(A,\O)\}$. Moreover, $\Pi(\zeta_A(\,\cdot\,,\O))=\{\re s>\ov\dim_B(A,\O)\}$; that is,
\begin{equation}
D(\zeta_A(\,\cdot\,,\O))=\ov\dim_B(A,\O).
\end{equation}

\medskip

$(b)$ If the box $($or Minkowski$)$ dimension $D:=\dim_B(A,\O)$ exists, $D<N$,  and $\M_*^D(A,\O)>0$, then $\zeta_A(s,\O)\to+\ty$ as $s\in\eR$
converges to $D$ from the right. In particular, $\mathcal{H}(\zeta_A(\,\cdot\,,\O))=\Pi(\zeta_A(\,\cdot\,,\O))=\{\re s>\dim_B(A,\O)\}$; that is,
\begin{equation}
D_{\rm hol}(\zeta_A(\,\cdot\,,\O))=D(\zeta_A(\,\cdot\,,\O))=\dim_B(A,\O).
\end{equation}
\end{theorem}

An entirely analogous result holds for the tube zeta function $\widetilde\zeta_A(\,\cdot\,,\O)$, except for the fact that the hypothesis $D<N$ is no longer needed in  the counterpart of part $(b)$ of Theorem \ref{an2}.

A very useful property of relative distance zeta functions is the following {\em scaling property}: for any RFD $(A,\O)$ and for any positive real number $\g$, we have
\begin{equation}\label{scaling}
\zeta_{\g A}(s,\g\O)=\g^s\zeta_A(s,\O).
\end{equation}

We refer the interested reader to \cite[Chap.\ 4]{fzf} and [LapRa\v Zu4--7]  for many other related results, examples and comments.
We mention, in particular, that `fractal drums' (that is, `drums with fractal boundary', in the sense of [Lap1--3], for example)\footnote{See also [Lap-vFr3, {\S}12.5], [Lap6] and [LapRa\v Zu7] for many other references on fractal drums.} correspond to RFDs of the form $(\pa\O,\O)$, where $\O$ is a nonempty bounded open subset of $\eR^N$, and that the results discussed in {\S}\ref{maxhf} above are applied in a crucial way in order to show the optimality of certain inequalities pertaining to the meromorphic continuations of the spectral zeta functions of fractal drums (viewed as RFDs); see \cite[{\S}4.3]{fzf} and \cite{brezish}.

\section{Fractal tube formulas and a Minkowski measurability criterion}\label{mmc}

In this section, we briefly explain how under suitable growth conditions on the relative distance (or tube) zeta function (see a variant of the languidity
(resp., of the strong languidity) condition of \cite[{\S}5.3]{lapidusfrank12} given in
~[LapRa\v Zu5,6]), it is possible to recover a pointwise or distributional fractal tube formula for a relative fractal drum $(A,\O)$ in $\eR^N$, expressed as a sum of residues over its visible complex dimensions.
These fractal tube formulas, along with a Tauberian theorem due to Wiener and Pitt (which generalizes Ikehara's Tauberian theorem, see~\cite{Kor, Pos}) make it possible to derive a Minkowski measurability criterion for a large class of relative fractal drums (and compact subsets) of $\eR^N$.
These results generalize to higher dimensions the corresponding ones obtained for fractal strings (that is, when $N=1$) in \cite[{\S}8.1 and {\S}8.3]{lapidusfrank12}.

The results of this section are announced in~\cite{cras1} and fully proved in~\cite{cras2}.
(See also \cite[Chap.\ 5]{fzf}.)
Furthermore, we refer the interested reader to [LapRa\v Zu1,6] and \cite[{\S}8.2 and {\S}13.1]{lapidusfrank12} for additional references on tube formulas in various settings, including [DeK\"O\"U, Fed1, Gra, HuLaWei, Schn, Z\"ah, Wey, LapPe1--2, LapPeWi1--2, LapLu, LapLu-vFr1--2]. (See also \cite[\S13.1, \S13.2 and \S13.4]{lapidusfrank12}.)

In order to be able to state the fractal tube formulas, we introduce the following notions, adapted from \cite{lapidusfrank12} to the present much more general context.
The {\em screen} $S$ is the graph of a bounded, real-valued, Lipschitz continuous function $S(\tau)$, with the horizontal and vertical axes interchanged:
$S:=\{S(\tau)+\I \tau\,:\,\tau\in\eR\}$
and we let $\sup S:=\sup_{\tau\in\eR}S(\tau)\in\eR$.
Given a relative fractal drum $(A,\O)$ of $\eR^N$, we always assume that the screen $S$ lies to the left of the critical line $\{\re s=\ov{\dim}_B(A,\O)\}$, i.e., that $\sup S\leq\ov{\dim}_B(A,\O)$.
Furthermore, the {\em window} $W$ is defined as $W:=\{s\in\Ce:\re s\geq S(\im s)\}$.
The relative fractal drum $(A,\O)$ is said to be {\em admissible} if its tube (or distance) zeta function can be meromorphically extended to an open connected neighborhood of some window $W$.

Assume now that $(A,\O)$ is an admissible relative fractal drum of $\eR^N$ for some screen $S$ such that its distance zeta function satisfies appropriate growth conditions (see [LapRa\v Zu5,6] for details).\footnote{Roughly speaking, $\zeta_{(A,\O)}:=\zeta_A(\,\cdot\,,\O)$ is assumed to grow at most polynomially along the vertical direction of the screen and along suitable horizontal directions (avoiding the poles of $\zeta_{(A,\O)}$); see \cite[Def.\ 5.2]{lapidusfrank12} for the so-called ``languidity condition''.}
Then its relative tube function satisfies the following identity, for all positive real numbers $t$ sufficiently small:\footnote{The ranges within which the formulas are valid are fully specified in [LapRa\v Zu5,6].}
\begin{equation}\label{distr_tube_A}
|A_t\cap\O|=\sum_{\omega\in\po(\zeta_A(\,\cdot\,,\O),W)}\res\left(\frac{t^{N-s}}{N-s}{\zeta}_A(s,\O),\omega\right)+\mathcal R(t).
\end{equation}

The above {\em fractal tube formula} is interpreted pointwise or distributionally, depending on the growth properties of $\zeta_A(\,\cdot\,,\O)$ and then, $\mathcal R(t)$ is a pointwise or distributional\footnote{For the precise definition of distributional asymptotics, see
\cite[{\S}5.4.2]{lapidusfrank12},  [LapRa\v Zu5--6]
and the relevant references therein.} asymptotic error term of order at most $O(t^{N-\sup S})$ as $t\to0^+$.
Moreover, if $S$ lies strictly to the left of the vertical line $\{\re s =\sup S\}$ (that is, if $S(\tau)<\sup S$ for every $\tau\in\eR$), then $\mathcal R(t)$ is $o(t^{N-\sup S})$, pointwise or distributionally, as $t\to 0^+$.   
In the case when $\zeta_A(\,\cdot\,,\O)$ satisfies stronger growth assumptions (i.e., the analog of the ``strong languidity conditon'' of \cite[Def.\ 5.3]{lapidusfrank12}), we obtain a tube formula without an error term (i.e., $\mathcal R(t)\equiv0$) and with $W=\Ce$.
Following [Lap-vFr3], the resulting formula is then called an {\em exact} fractal tube formula.

The tube formula \eqref{distr_tube_A} can also be expressed in terms of the relative tube zeta function when analogous growth conditions are imposed on $\widetilde{\zeta}_A(\,\cdot\,,\O)$:\footnote{Note that in light of the functional equation \eqref{equ_tilde}, assuming growth conditions for $\zeta_A$ is essentially equivalent to assuming them for $\widetilde\zeta_A$ (and vice versa).}
\begin{equation}\label{distr_tube_tube_A}
|A_t\cap\O|=\sum_{\omega\in\po(\widetilde{\zeta}_A(\,\cdot\,,\O),W)}\res\left({t^{N-s}}\widetilde{\zeta}_A(s,\O),\omega\right)+\mathcal R(t).
\end{equation}
In fact, the key observation for deriving the above formula is the fact that
\begin{equation}
\widetilde{\zeta}_A(s,\O)=\int_0^{+\ty}t^{s-N-1}\chi_{(0,\d)}(t)|A_t\cap\O|\,\D t=\{\mathfrak{M}f\}(s),
\end{equation}
where $\chi_E$ is the characteristic function of the set $E$, $\{\mathfrak{M}\psi\}(s):=\int_0^{+\ty}t^{s-1}\psi(t)\,\D t$ is the Mellin transform of the function $\psi$, and $f(t):=t^{-N}\chi_{(0,\d)}(t)|A_t\cap\O|$.
One then applies the inverse Mellin transform (see \cite{titch}) to recover the relative tube function $t\mapsto|A_t\cap\O|$ and  proceeds in a similar manner as in~\cite[Chap.\ 5]{lapidusfrank12} for the case of fractal strings.

As an application, the following result generalizes the Minkowski measurability criterion given in \cite[Thm.\ 8.15]{lapidusfrank12} for fractal strings to the present case of relative fractal drums.

\begin{theorem}[Minkowski measurability criterion, {\rm [LapRa\v Zu5,6]}]\label{criterion}
Let $(A,\O)$ be an admissible relative fractal drum of $\eR^N$ such that $D:=\dim_BA$ exists and $D<N$.
Furthermore, assume that its relative distance $($or tube$)$ zeta function satisfies appropriate growth conditions\footnote{See [LapRa\v Zu5,6] for details about these growth conditions.} for a screen passing between the critical line $\{\re s=D\}$ and all the complex dimensions of $A$ with real part strictly less than $D$.
Then, the following statements are equivalent$:$

\medskip

$(a)$ $A$ is Minkowski measurable.

\medskip

$(b)$ $D$ is the only pole of the distance zeta function ${\zeta}_A$ located on the critical line $\{\re s=D\}$, and it is simple.
\end{theorem}

There exist relative fractal drums which do not satisfy the hypothesis of Theorem~\ref{criterion} concerning the screen; see [Lap-vFr3, Exple.\ 5.32]. 
We point out that the fractal tube formula \eqref{distr_tube_A} can be used to recover (or obtain for the first time) the (relative) fractal tube formulas for a variety of well-known (and not necessarily self-similar) fractal sets, as is illustrated by the following examples.

\begin{example}
Recall from Proposition~\ref{sierpinski_carpet0} that the distance zeta function of the Sierpi\'nski carpet $A$ is given for all $s\in\Ce$ by
$$
\zeta_A(s)=\frac{8}{2^ss(s-1)(3^s-8)}+2\pi\frac{\d^s}s+4\frac{\d^{s-1}}{s-1},
$$
for $\d>1/6$, and is meromorphic on all of $\Ce$.
It is easy to check that $\zeta_A$ satisfies growth conditions which are good enough for \eqref{distr_tube_A} to hold pointwise without an error term and for all $t\in(0,1/2)$:
\begin{equation}\label{carpp}
|A_t|=\sum_{{\omega\in\po(\zeta_A,\Ce)}}\res\left(\frac{t^{2-s}}{2-s}{\zeta}_A(s),\omega\right).
\end{equation}
Now, also recall from Proposition \ref{sierpinski_carpet0} that $\po(\zeta_A,\Ce)=\{0,1\}\cup\{s_k:k\in\Ze\}$, where $s_k=\log_38+\frac{2\pi}{\log 3}k\I$ for all $k\in\Ze$.
Furthermore, $\res(\zeta_A,0)=2\pi+8/7$, $\res(\zeta_A,1)=16/5$ and the residues at $s_k$ are given in \eqref{ressk}; so that \eqref{carpp} becomes the following exact, pointwise fractal tube formula, valid for all $t\in(0,1/2)$:
\begin{equation}
|A_t|=\frac{t^{2-\log_38}}{\log 3}\sum_{k=-\ty}^{+\ty}\frac{2^{-s_k}t^{-\frac{2\pi}{\log 3}k\I}}{s_k(s_k-1)(2-s_k)}+\frac{16}{5}t+\left(2\pi+\frac{8}{7}\right)t^2.
\end{equation}
\end{example}

The above example can be generalized to an $N$-dimensional analog of the Sierpi\'nski carpet (see [LapRa\v Zu1,6]).
We next establish the special case of this assertion for the relative $3$-dimensional Sierpi\'nski carpet.

\begin{example}\label{ex2}
Let $A$ be the three-dimensional analog of the Sierpi\'nski carpet and $\O$ the closed unit cube in $\eR^3$.
More precisely, we construct $A$ by dividing $\O$ into $27$ congruent cubes and remove the open middle cube, then we iterate this step with each of the $26$ remaining smaller closed cubes; and so on, ad infinitum.
By choosing $\d>1/6$, we deduce that $\zeta_A$ is meromorphic on $\Ce$ and given for all $s\in\Ce$ by (see~[LapRa\v Zu1,6])
\begin{equation}
\zeta_A(s,\O)=\frac{48\cdot 2^{-s}}{s(s-1)(s-2)(3^s-26)}.
\end{equation}
In particular, $\po({\zeta}_A(\,\cdot\,,\O),\Ce)=\{0,1,2\}\cup\big(\log_326+\mathbf{p}\I\Ze\big)$, where $\mathbf{p}:=2\pi/\log 3$.
Furthermore, we have that
$$
\res(\zeta_A(\,\cdot\,,\O),0)=-\frac{24}{25},\q \res(\zeta_A(\,\cdot\,\O),1)=\frac{24}{23},\q \res(\zeta_A(\,\cdot\,,\O),2)=-\frac{6}{17}
$$
and, by letting $\omega_k:=\log_326+\mathbf{p}k\I$ for all $k\in\Ze$, 
$$
\res(\zeta_A(\,\cdot\,,\O),\omega_k)=\frac{24\cdot 2^{-\omega_k}}{13\cdot \omega_k(\omega_k-1)(\omega_k-2)\log 3}.
$$
Again, the relative distance zeta function $\zeta_A(\,\cdot\,,\O)$ satisfies sufficiently good growth conditions, which enables us to obtain the following exact pointwise relative tube formula, valid for all $t\in(0,1/2)$:
$$
\begin{aligned}
|A_t\cap\O|&=\frac{24\,t^{3-\log_326}}{13\log 3}\sum_{k=-\ty}^{+\ty}\frac{2^{-\omega_k}t^{-\mathbf{p}k\I}}{(3-\omega_k)(\omega_k-1)(\omega_k-2)\omega_k}-\frac{6}{17}t+\frac{12}{23}t^2-\frac{8}{25}t^3.
\end{aligned}
$$
In particular, we conclude that $\dim_B(A,\O)=\log_326$ and, by Theorem \ref{criterion}, that, as expected, $(A,\O)$ is not Minkowski measurable.
\end{example}

One can similarly recover the well-known fractal tube formula for the Sierpi\'n\-ski gasket obtained in~\cite{lappe2} (and also, more recently, by a somewhat different method, in~\cite{DeKoOzUr}), as well as a tube formula for its $N$-dimensional analog described in \cite[Chap.\ 5]{fzf}.\footnote{We can also recover and extend the significantly more general fractal tube formulas obtained (for fractal sprays and self-similar tilings) in \cite{lappewi1} and used, in particular, in \cite{lappewi2}.}
We also point out that, in light of the functional equation \eqref{dist_geo}, the above fractal tube formulas \eqref{distr_tube_A} and \eqref{distr_tube_tube_A} generalize the corresponding ones obtained for fractal strings (i.e., when $N=1$) in \cite[{\S}8.1]{lapidusfrank12}.
Furthermore, these tube formulas can also be applied to a variety of fractal sets that are not self-similar, including `fractal nests' and `geometric chirps' (see \cite[Chaps.\ 3 and 4]{fzf} for the definitions of these notions and \cite{cras2} along with \cite[Chap.\ 5]{fzf} for the actual fractal tube formulas).

We conclude this section by briefly explaining how these results can also be applied in order to recover (and extend) the tube formulas for {\em self-similar sprays} generated by a suitable bounded open set $G\subset\eR^N$.
(See \cite{lappe2, lappewi1}.)
A self-similar spray is a collection $(G_k)_{k\in\eN}$ of pairwise disjoint sets $G_k\subset\eR^N$, with $G_0:=G$ and such that $G_k$ is a scaled copy of $G$ by some factor $\lambda_k>0$.
The sequence $(\lambda_k)_{k\in\eN}$ is called the {\em scaling sequence} associated with the spray  and is obtained from a ``ratio list'' $\{r_1,r_2,\ldots,r_J\}$, with $0<r_j<1$ for each $j\in\{1,2,\dots,J\}$, by building all possible words based on the ratios $r_j$.
Let now $(A,\O)$ be the relative fractal drum such that $A:=\partial(\cup_{k\in\eN} G_k)$ and $\O:=\cup_{k\in\eN} G_k$, with $\ov{\dim}_B(\partial G,G)<N$.
Then, it is clear that its relative distance zeta function $\zeta_A(\,\cdot\,,\O)$ satisfies the following functional equation, for all $s\in\Ce$ with $\re s$ sufficiently large:
\begin{equation}
\zeta_A(s,\O)=\zeta_{\partial G}(s,G)+\zeta_{r_1A}(s,r_1\O)+\cdots+\zeta_{r_JA}(s,r_J\O),
\end{equation}
where $(r_jA,r_j\O)$ denotes the relative fractal drum $(A,\O)$ scaled by the factor $r_j$.
Furthermore, by using the scaling property \eqref{scaling} of the relative distance zeta function, the above equation becomes
\begin{equation}
\zeta_A(s,\O)=\zeta_{\partial G}(s,G)+r_1^s\zeta_{A}(s,\O)+\cdots+r_J^s\zeta_{A}(s,\O),
\end{equation}
which yields that
\begin{equation}
\zeta_A(s,\O)=\frac{\zeta_{\partial G}(s,G)}{1-\sum_{j=1}^{J}r_j^s}.
\end{equation}
It is now enough to assume that the relative distance zeta function $\zeta_{\partial G}(s,G)$ of the generating relative fractal drum $(\partial G,G)$ satisfies suitable growth conditions in order to obtain the following formula for the `inner' volume of $A_t=(\pa\O)_t$ relative to $\O:=\cup_{k\in\eN} G_k$, for all positive $t$ sufficiently small:
\begin{equation}\label{exact}
|A_t\cap\O|=\sum_{\omega\in\mathfrak{D}(W)\cup\po(\zeta_{\partial G}(\,\cdot\,,G),W)}\res\left(\frac{t^{N-s}{\zeta}_{\partial G}(s,G)}{(N-s)\Big(1-\sum_{j=1}^Jr_j^s\Big)},\omega\right)+\mathcal R(t),
\end{equation}
where $\mathfrak{D}(W)$ denotes the set of all visible complex solutions of $\sum_{j=1}^Jr_j^s=1$ (in $W$) and $W$ is the window defined earlier.
It is easy to check that (at least) in the case of monophase or pluriphase generators (in the sense of \cite{lappe2} and [LapPeWi1,2]), these growth conditions are satisfied, so that one obtains exactly the same distributional or pointwise fractal tube formulas as in \cite{lappe2} or \cite{lappewi1}, respectively, after having calculated the distance zeta function $\zeta_{\pa G}(\,\cdot\,,G)$ of the generator. Moreover, if $\zeta_{\pa G}(\,\cdot\,,G)$ is strongly languid, we can let ${\mathcal R}(t)\equiv0$ and $W=\Ce$ in \eqref{exact} and therefore obtain {\em exact} fractal tube formulas.
\smallskip

We conclude this survey by pointing out that a broad variety of open problems and suggestions for directions for future research in this area are proposed in \cite[Chap.\ 6]{fzf}.

\end{document}